\documentclass{amsart}

\usepackage{epsfig}
\usepackage{latexsym}
\usepackage{amsmath}

\newtheorem{theorem}{Theorem}[section]
\newtheorem{proposition}[theorem]{Proposition}

\newtheorem{lemma}[theorem]{Lemma}
\newtheorem{corollary}[theorem]{Corollary}

\newcommand{\RR}{{\mathbb R}}
\newcommand{\PP}{{\mathbb P}}
\newcommand{\EE}{{\mathbb E}}
\newcommand{\II}{{\mathbb I}}

\newcommand{\old}[1]{{}}

\title{Inverting random functions III: 
discrete MLE revisited}
\author{Mike A. Steel and L\'aszl\'o A. Sz{\'e}kely}
\thanks{We thank the NZIMA (Maclaurin
Fellowship) for supporting this research. The second author was
also supported in part by NSF DMS contract 007 2187.}

\address{Biomathematics Research Centre, Department of Mathematics and
  Statistics, University of Canterbury, Christchurch, New Zealand. \\
  Department of Mathematics, University of South
   Carolina \\ Columbia SC, USA.}

\email{m.steel@math.canterbury.ac.nz, szekely@math.sc.edu}


\date{3 August 2006}


\begin{document}
\begin{abstract}
This paper continues our earlier investigations into the inversion of random
functions in a
general (abstract) setting. In Section~\ref{two} we investigate a concept of
invertibility and the invertibility of the composition of random
functions. In  Section~\ref{three} we
resolve some questions concerning the number of samples required to ensure the
accuracy of
parametric maximum likelihood estimation (MLE). A direct application to 
phylogeny reconstruction is given.
\end{abstract}

\maketitle

\section{Review of random functions}
This paper is a sequel of our earlier papers \cite{inv1, inv2}. We assume that
the reader is familiar with those papers; however, we repeat the most important
definitions.

For two finite
sets, $A$ and $U$, let us be given a $U$-valued random variable
$\xi_a$ for every $a\in A$. We call the vector of random variables
$(\xi_a: \ a\in A)$ a {\em random function}
$\Xi:\  A\rightarrow U$. 
Ordinary functions are specific
instances of random functions.

Given another random function, $\Gamma$, from $U$ to $V$, we can speak about
the composition of $\Gamma$ and $\Xi$, $\Gamma\circ\Xi: \ A\rightarrow V$,
 which is the vector variable 
$(\gamma_{\xi_a}: \ a\in A)$.
In this paper we are concerned with inverting 
random functions. In other words, we look for random functions
$\Gamma: \ U \rightarrow A$ in order to obtain the best approximations of the 
identity function $\iota:\ A\rightarrow A$ by  $\Gamma\circ\Xi$. 
{\em We always assume that
$\Xi$ and $\Gamma $ are independent.} This assumption holds for free if
either $\Xi$ or $\Gamma $ is a deterministic function.

Consider the probability of returning $a$ from $a$ by the
composition of two random functions, that is,
$r_a=\PP[\gamma_{\xi_a}=a]$. The assumption on the independence of
$\Xi$ and $\Gamma$ immediately implies
\begin{equation}\label{indeps} r_a=\sum_{u\in U} \PP[\xi_a=u]\cdot
\PP[\gamma_u=a].  \end{equation} A natural criterion is to find
$\Gamma$ for a given $\Xi$ in order to maximize $\sum_a r_a$.  More
generally, we may have a weight function
$w:A\rightarrow \RR^+$ and we may wish to maximize $\sum_a r_aw(a)$. This
can happen if we give preference to returning certain $a$'s, or, if
we have a  prior probability distribution on $A$ and we want to
maximize the expected return probability for a random element of $A$
selected according to the prior distribution.  The following random function 
$\Gamma^*:U\rightarrow A$, defined below, will do this job:
 for any fixed $u\in U$,
\begin{equation}\label{mlenonpar} \gamma_u^*=a^* \hbox{\ \rm for
sure, if for all\ } a\in A, \ \ \PP[\xi_{a^*}=u]w(a^*)\geq
\PP[\xi_{a}=u]w(a).  \end{equation} 
(In case there is more than one element $a^*$
that satisfies (\ref{mlenonpar}), we may select uniformly at random
from the set of such elements.)
This function $\Gamma^*$  is
called  the {\em maximum a posteriori estimator} (MAP) in the
literature \cite{everitt98}. 
The special case when the weight function $w$ is constant, is known as
the {\em maximum likelihood estimation} (MLE) \cite{casellaberger,
  everitt98}.

For $a,b\in A$, $\Xi:A\rightarrow U$, let
\begin{equation} \label{defvardist}
d(a,b)=: d(\xi_a,\xi_b)=   \sum_{u\in U} \biggl\vert \PP[\xi_a=u]-\PP[\xi_b=u]
\biggl\vert,
\end{equation}
which is called the {\em variational distance} of 
the random variables
$\xi_a$ and $\xi_b$. 

A given $\Xi: \ A\rightarrow U$  will have an $|A|\times |U|$
{\em associated matrix} $X$, such that $x_{au}=\PP[\xi_a=u]$.
Given a $\Gamma :  \  U\rightarrow V$ with associated matrix $G$, 
the composition of $\Gamma$ and $\Xi$, $\Gamma\circ\Xi: \ A\rightarrow V$,
will have the associated matrix $XG^T$.

Our motivation for the study of random functions came from phylogeny
reconstruction \cite{fels, sem03}. Stochastic models
 define how biomolecular sequences 
are generated at the leaves of a binary tree.  If all possible binary
trees on $n$ leaves come equipped with a model for generating
biomolecular sequences of length $k$, then we have a random function
from the set of binary trees with $n$ leaves to the ordered $n$-tuples
of biomolecular sequences of length $k$. {\em Phylogeny reconstruction}
can be viewed as a random function from the  set of ordered $n$-tuples
of biomolecular sequences of length $k$ to the  set of binary trees
 with $n$ leaves. It is a natural assumption that random mutations in the
past are independent from any random choices in the  phylogeny
reconstruction algorithm. Criteria for phylogeny
reconstruction may differ according to what one wishes to optimize.
However, in the practice of phylogeny reconstruction there are no fixed,
preconceived models on the possible trees; instead, we also try to find 
out the model parameters.
Our paper \cite{inv1} introduced a new abstract model for phylogeny
reconstruction: inverting parametric random functions.  Most of the work
done on the mathematics of phylogeny reconstruction can be discussed in
this context.  This model is more structured than random functions, and
hence is better suited to describe details of models of phylogeny and
the evolution of biomolecular sequences.

Assume that for a finite set $A$, for every $a\in A$, an 
(arbitrary, finite or infinite) set
$\Theta(a)\not=\emptyset $ is  assigned, and moreover, $\Theta(a)\cap
\Theta(b)=\emptyset$  for $a\not= b$. Set 
$B=\{  (a, \theta) :\ a\in A, \theta\in \Theta(a)\}$ and let $\pi_1$
 denote the 
natural
projection from $B$ to $A$.  A {\em parametric random function} is
the collection $\Xi$ of random variables such that

\noindent  for $a\in A$ and $\theta\in \Theta(a)$, there is a (unique)
$U$-valued  random variable $\xi_{(a, \theta)}$  in $\Xi$.

\begin{figure}[h]
\begin{center} \label{hybridfig}
\resizebox{12cm}{!}{
\input{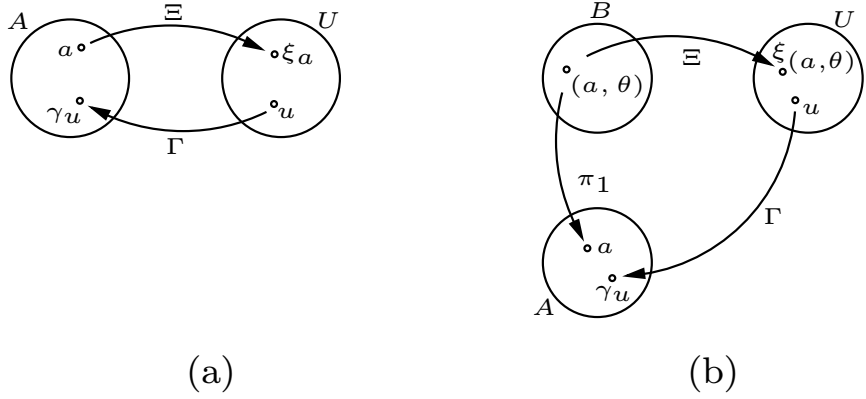}
}
\caption{Inversion of non-parametric (a), and parametric (b) random functions}
\end{center}
\label{overview}
\end{figure}

We are interested in random functions $\Gamma:U \rightarrow A$
independent from $\Xi$ so that $\gamma_{\xi_{(a,\theta)}}$ best
approximates $\pi_1$ under certain criteria. Call $R_{(a,\theta)}$
the probability $\PP[\gamma_{\xi_{(a,\theta)}}=a]$. Maximum
Likelihood Estimation, as it is used in situations where there  is
 a discrete parameter of interest to estimate,
 in the presence of other parameters (such as phylogeny
reconstruction), would take the $\Gamma'$, for which for every fixed
$u$, $\gamma_u'=a'$ for sure, if \begin{equation} \label{mlepara} \
\forall (a,\theta)\in B\ \ \exists \theta' \in \Theta(a') \ \
\PP[\xi_{(a',\theta')}=u]\geq \PP[\xi_{(a,\theta)}=u]. 
\end{equation}
In case there is more than one element $a'$ that satisfies
(\ref{mlepara}), we may select uniformly at random from the set of
such elements. (We avoided using the more natural looking
quantification  $ \exists \theta' \in
\Theta(a') \ \ \forall (a,\theta)\in B$, since
$\PP[\xi_{(a',\theta')}=u]$ may not take a maximum value!) We denote
by $R'_{(a, \theta)}$ the probability that from the pair
$(a,\theta)$ the  Maximum Likelihood Estimation $\Gamma'$ returns
$a$, i.e.  \begin{equation} \label{mleR}
R'_{(a,\theta)}=\PP[\gamma'_{\xi_{(a,\theta)}}=a].  \end{equation}

If a random function 
$\Xi:\  A\rightarrow U$ ($\Xi:\  B\rightarrow U$)
  is to have $k$ independent evaluation, we denote
the resulting random function by $\Xi^{(k)}:\  A\rightarrow U^k$
($\Xi^{(k)}:\  B\rightarrow U^k$), and  the
random variable associated with $a$ will be $\xi^{(k)}_a$. We will
study the invertibility of $\Xi^{(k)}$ both in the non-parametric and the 
parametric setting. For a $\Gamma: U^k\rightarrow A$ random function, we 
use the notation $r^{(k)}_a=\PP[\gamma_{\xi_a^{(k)}}=a]$ in the non-parametric
case, $R^{(k)}_{(a,\theta)}=\PP[\gamma_{\xi_{(a,\theta)}^{(k)}}=a]$ in the
parametric case, and $[R^{(k)}]'_{(a,\theta)}$, if $\Gamma'$ is the Maximum
Likelihood Estimation.

In Section~\ref{two} we will show that in the non-parametric setting
several natural definitions of invertibility of a random function are,
in fact, equivalent. Furthermore, we determine  when composition of 
invertible random functions is invertible. The main result of this Section
is an explicit bound on how invertibility ``improves'' as the variational
distances between elements of $A$ have increasing separation from zero.

In Section~\ref{three} we revisit our study of
the worst-case behavior of MLE in \cite{inv2}. (This is a
very natural question in situations where a prior distribution is
not given on $A$, or the inverting of the random function is to be
carried out only once. Such a situation arises in phylogeny
reconstruction, where, arguably, we do not have a prior distribution on
alternative evolutionary scenarios, and the reconstruction is not
going to be repeated---there is only one `Tree of Life' that we want
to know.)   A certain amount of controversy and debate has surrounded the
statistical consistency of MLE in phylogeny, as described in 
\cite{fels}, pp. 270--272. Felsenstein's claim (from the early 1970s) of the
 consistency of MLE in 
phylogeny for simple (`identifyable') models is correct, but it
 was only formally established in 1996 by 
\cite{chang}. This result, like
Wald's earlier result \cite{wald}, relies on a compactness
argument, continuity, and limit theory,
that does not give an explicit bound on $k$. Other proofs in the
 biological literature
have generally been less rigorous and led to criticism and debate
 (see eg. \cite
{far99, rog97, rog01, sid98, yan94, yan96}). One oversight has been
 to treat the 
MLE-estimated continuous parameters (branch lengths) of alternative
 trees as fix
ed rather than as random variables dependent on the data; such 
arguments are satisfying for practical purposes but call for more rigor.
The significance of Theorem 5.1 \cite{inv2}  is that it gives the first
 explicit bounds for MLE, both in the phylogenetic setting and
beyond.  However, this result depended on an unnatural parameter, namely
the smallest positive probability  that an image of the object to be 
reconstructed can have. Here in Theorem~\ref{mainpara}
 we get rid of this dependence, and provide
a simple and immediate application of this new result to phylogeny
 reconstruction.

We study two examples that show how subtle is MLE for  inverting
parametric random
functions. The first example shows that Theorem~\ref{mainpara} is
``near optimal'' in one of its parameters.   The second example shows that 
in contrast to the non-parametric setting, the vanishing of
 variational distance
does not by itself preclude MLE (or other) estimation for certain 
random functions.

Our approach is information-theoretic, we focus on the possibility or
impossibility of inverting random functions, and not on the computational
complexity issues. Our results can also be re-stated in the language of
decision theory, by talking about `loss functions' and `risk function'
associated to the  decision rule.

\section{invertibility in the non-parametric setting} \label{two}
Let us say that a random function $\Xi:A \rightarrow U$  is {\em invertible}
 if there exists a random function
$\Gamma: U \rightarrow A$  such that
for all $a \in A$,
$\PP[\gamma_{\xi_a} = x]$ takes strict maximum when $x=a$, or
equivalently,
\begin{equation}
\label{difference}
\PP[\gamma_{\xi_a} = a] - \max_{x \neq a} \{\PP[\gamma_{\xi_a} =x]\} >0
\mbox{ for all  $a \in A$}.
\end{equation}
Informally, $\Xi$ is
invertible, if there is some reconstruction method that is always more
likely to pick the generating object in $A$ than any other
element of $A$.

A sufficient condition for $\Xi$ to be  invertible is that there
exists a $\Gamma$ so that for all $a \in A$,  the following two
conditions apply:
\begin{itemize}
\item[($I_1$)]
$\PP[\gamma_{\xi_a} = a] > \frac{1}{|A|},$
\item[($I_2$)] $\PP[\gamma_{\xi_a} = b] < \frac{1}{|A|}, \mbox{ for all
} b \neq a.$
\end{itemize}

Note that invertibility implies ($I_1$), and is equivalent to it 
when $|A|=2$, but not equivalent for $|A|\geq 3$.

We say $\Xi$ {\em separates} $A$, if, for each distinct pair
$a,b \in A$,
the variational distance $d(a,b)$ of the probability distributions
of
$\xi_a$ and $\xi_b$ is strictly positive.

\begin{proposition}
\label{main1}
\mbox{}
The following properties are equivalent for an $\Xi:\ A\rightarrow U$ 
random function:
\begin{itemize}
\item[(i)]
$\Xi$ separates $A$
\item[(ii)]
For all $\epsilon>0$
there is a value of $k_{\epsilon}$ so that for all
$k \geq k_{\epsilon}$ there is a random function $\Gamma^\S: \ U^k
\rightarrow A$ for which
$\PP[\gamma_{\xi^{(k)}_a}^\S=a] > 1-\epsilon$.
\item[(iii)]
$\Xi$ is  invertible
\item[(iv)]
For some $k\geq 1$, $\Xi^{(k)}$ is  invertible.
\end{itemize}
\end{proposition}
\begin{proof}
The equivalence between (i) and (ii) follows easily from results in our 
earlier
papers \cite{inv1}
and \cite{inv2} and standard arguments.  We will show that
(iv) $\Rightarrow$ (ii) and that (i) $\Rightarrow$ (iii). Since (iii)
$\Rightarrow$ (iv) is trivial this will establish the claimed four-way
equivalence.

{\em Proof of (iv) $\Rightarrow$ (ii)}
Suppose that $\Xi^{(k)}$ is  invertible.
Select $\Gamma$ to satisfy (\ref{difference}) for $\Xi^{(k)}$.
For positive integer $m$, generate $km$ independent samples in $U$ 
according to
$\Xi$. Define
$\Gamma^\S :\ U^k\rightarrow A$ as follows: select the elements of $A$ that are
reconstructed
most often according to $\Gamma$ and choose one of them uniformly at
random. By standard probability arguments, the
probability that the correct element $a$ will be selected by this
process converges to 1 as $m$ tends to infinity.

{\em Proof of (i) $\Rightarrow (iii)$}
Suppose that $\Xi:\ A\rightarrow U$ separates $A$.
Let $X$ denote the associated matrix of  $\Xi$,
and let ${\mathbf a}_i$, $i\in A$ denote the rows of $X$. Recall that
${\mathbf a}_i$ gives the distribution of $\xi_i$.
We will describe the inverse random   function $\Gamma:\ U\rightarrow A$
with its associated matrix, i.e.
in the form of a $|U|\times |A|$ matrix $G$, whose rows represent
the distribution of the element of $U$ corresponding to the row.

We write  $G=V+{1\over |A|} J$ and will give $V$ explicitly.
(If we
were to
take $V=0$, then (\ref{difference}) yields uniformly $=0$ instead of 
the desired
$>0$). We denote the columns of $V$ by  ${\mathbf v}_i$, $i\in A$.
We define each vector ${\bf v}_i$ as follows:
  $${\bf v}_i=  \frac{{\bf a}_i}{|{\bf a}_i|}-  \frac{1}{|A|} \sum_{j=1}^{|A|}
\frac{{\bf a}_j}{|{\bf a}_j|},$$
where $|.|$ is the usual euclidean vector norm. 
 Then it can be checked that this
choice of $V$ provides a solution to
the following system:
\begin{eqnarray*}
\forall i \forall j\not= i\ \ \  {\mathbf a}_i\cdot  {\mathbf v}_i 
-{\mathbf a}_i\cdot  {\mathbf v}_j-\epsilon_{ij} & = & 0; \\
\sum_{l\in A}  {\mathbf v}_l & = & 0; \\ 
\forall i \forall j\not= i\ \ \  \epsilon_{ij} & > & 0. 
\end{eqnarray*}
and these are precisely the conditions (\ref{difference}) requires
 for  invertibility.
\end{proof}

\subsection{Composition of invertible functions}

A natural question is whether the composition of 
 invertible functions is also  invertible.
The next result shows that in general the answer is `no', though 
we can provide
  a precise characterization based on the
rank of an associated matrix.

\begin{theorem}
Let  $\Upsilon: U\rightarrow Z$ be a random function 
matrix $Y$, and let $Y^+$  denote the extension of $Y$ by an all-1 row.
If  $rank(Y^+)=|U|$,  then for all
   $\Xi: A\rightarrow U$
 invertible
random functions, the composition $\Upsilon \circ \Xi: A\rightarrow Z$
is  invertible, and if   rank is less
than $|U|$, then there exist  invertible
random functions $\Xi: A\rightarrow U$ such that
$\Upsilon \circ \Xi: A\rightarrow Z$ is not  invertible.
\end{theorem}

\begin{proof}
Assume first that $\Upsilon \circ \Xi$ is not 
invertible, i.e. there exist $a\not= b\in A$, such that the
distributions $\upsilon_{\xi_a}$ and $\upsilon_{\xi_b}$ are identical.
Then we have the following homogeneous system of linear
equations, where the coefficients are the numbers   $\PP[\upsilon_u=z]$
and
1's,
and the variables are the $x_u$'s:
\begin{eqnarray} \label{identity}
\sum_{u\in U} \PP[\upsilon_u=z]x_u &= &0
\hbox{\ \ \  for all $z\in Z$.}\\
\sum_{u\in U} x_u& = & 0. \label{identity2}
\end{eqnarray}
The matrix $Y^+$ is the matrix of the system of homogeneous linear 
equations (\ref{identity})-(\ref{identity2}).
Observe that $x_u= \PP[\xi_a=u]-\PP[\xi_b=u]$ solves the system
(\ref{identity})-(\ref{identity2}).
   If the rank of $Y^+$
is $|U|$, then it has only trivial solution, i.e. for all $u\in U$
$x_u=0$. This amounts to $\xi_a$ and $\xi_b$ having the same
distribution,
contrary to the assumption of $\Xi$ being  invertible.

Assume now that  $Y^+$ has rank less than $|U|$. Then
the system (\ref{identity})-(\ref{identity2})
   has a non-trivial solution $x_u$.
Set $P=\sum_{u: \ x_u>0} x_u$ and $N=\sum_{u: \ x_u<0} x_u$. Clearly
$P=-N>0$. Take $A=\{a,b\}$,
$ \PP[\xi_a=u]= {x_u\over P} $  if $ x_u\geq 0$, and
0  otherwise; and $ \PP[\xi_b=u]= {x_u\over N} $  if $ x_u\leq 0$, and
0  otherwise. It is clear that this $\Xi $ is  invertible, as it 
separates $a$ and $b$.
However, according to the argument above (\ref{identity}), the
distributions $\upsilon_{\xi_a}$ and $\upsilon_{\xi_b}$ are identical.
\end{proof}

\subsection{Explicit bounds}

 From Proposition~\ref{main1}, if $\Xi$ separates $A$ then
there is
a random function $\Gamma: U \rightarrow A$ for which
$$\PP[\gamma_{\xi_a} = a] -\frac{1}{|A|} > 0.$$
We now consider putting an explicit lower bound on the right hand side of this
inequality.
That is, we show that for a specific continuous positive function $h: \RR
\rightarrow \RR$
(dependent only on
$|A|$) the following holds:
Suppose that $d(a,b) > \delta$ for all $a, b \in A, a \neq b$.  Then
there is
a random function $\Gamma: U \rightarrow A$ for which
$$\PP[\gamma_{\xi_a} = a] -\frac{1}{|A|} > h(\delta)$$ for all $a \in A$.
Note that we cannot insist the $\Gamma$ be MLE 
(maximum likelihood estimation), even when
$|A|=2$. To see this, let
$A=\{1,2\}, U = \{u_1, u_2\}$ and let $\xi_1$ take the value $u_1$ with
probability 1, and let $\xi_2$ take the values $u_1, u_2$ with probabilities
$\frac{2}{3}$ and
$\frac{1}{3}$, respectively; then if $\Gamma=\Gamma^*$ is MLE, 
 we have $\PP[\gamma_{\xi_2}
= 2] = \frac{1}{3}$.

\begin{theorem}
For every random function
$\Xi:A \rightarrow U$, with $|A|>1$, there exists a $\Gamma: U\rightarrow A$,
such
that
\begin{equation}
\min_{a\in A} r_a \geq {1\over |A|} +{1\over 2|A|(|A|-1)}\min_{a\in A}
\sum_{b\in A} d(a,b).
\end{equation}
In particular, if for all $a\not= b \in A$, $d(a,b)\geq \delta$, then
$\min_{a\in A} r_a \geq {1\over |A|} +{\delta \over 2|A|}$. 
\end{theorem}
\begin{proof}
Recall the characterization of the random inverse function
maximizing $\min_{a\in A} r_a$ from Theorem 5 \cite{inv1}:
$\min_{a\in A} r_a= \min_\mu \sum_{u\in U} \max_{a\in A} \mu(a)
\PP[\xi_a=u]$,
where $\mu$ is a probability distribution on $A$. In the rest of the 
proof $\mu$ refers to this minimizing distribution.
(Note that Theorem 5 in
   \cite{inv1} contains an annoying typo, it shows maximization for $\mu$
instead of minimization). We are going to use the following Lemma.
\begin{lemma}
\label{bnums}
Let us be given real numbers $b_1,b_2,...,b_n$. Assume that
$$\sum_{1\leq i < j \leq n} |b_i-b_j|\geq (n-1)\epsilon.$$ Then
$\max_j [b_j -{1\over n} \sum_{i=1}^n b_i] \geq {\epsilon \over n}.$
\end{lemma}
\begin{proof}
Without loss of generality we may assume $b_1\geq b_2\geq ...\geq
b_n$. The conditions of the Lemma can be rewritten as the conditions
of the following primal linear program:
\begin{eqnarray*}
b_2-b_1 & \leq & 0 \\
b_3-b_2 & \leq & 0 \\
&...&\\
b_n-b_{n-1} & \leq & 0 \\
\sum_{i< j} b_i-b_j  & \leq & -(n-1)\epsilon\\
& & \max ({1\over n} \sum_i b_i) -b_1.
\end{eqnarray*}
Recall the Duality Theorem of linear programming \cite{schrijver}:
$\max\{c^Tx:\ Mx\leq b\}=\min\{y^Tb:\ y\geq 0, \ y^tM=c\},$
if both optimization problems have feasible solutions. The dual  linear
program  is as follows:
\begin{eqnarray*}
(n-1)x_n-x_1 & = & -{n-1\over n}\\
x_i-x_{i+1}+(n-2i-1)x_n& =&{1\over n} \ \ \hbox{for \ } i=1,2,...,n-2;\\
x_{n-1}+(1-n)x_n&=&  {1\over n}\\
x_1,x_2,...,x_n&\geq & 0\\
& & \min -(n-1)\epsilon x_n.
\end{eqnarray*}
It is easy to see that the for the dual problem a feasible solution  is
   the following
setting: $x_i=1-{i(i-1)\over n(n-1)}$ for $i=1,2,...,n-1$, and
    $x_n={1\over n(n-1)}$; with value $-{\epsilon\over n}$.
This implies that ${\epsilon\over n}\leq
   \max_j b_j -{1\over n} \sum_{i=1}^n b_i$ for any feasible solution of
   the primal problem.
\end{proof}
We are going to apply Lemma~\ref{bnums} in the following setting. Fix an
arbitrary $u\in U$, and  for $i\in A$,
let $b_i=\mu(i)\PP[\xi_i=u]$. The lemma yields
\begin{eqnarray} \label{nyolc}
& &\max_{a\in A} \Biggl(\mu(a)\PP[\xi_a=u]-{1\over |A|}\sum_{i\in A}
\mu(i)\PP[\xi_i=u]\Biggl) \\
& \geq & {1\over |A|(|A|-1)} \sum_{1\leq i < j \leq |A|}
\bigl|\mu(i)\PP[\xi_i=u]-\mu(j)\PP[\xi_j=u]\bigl|. \label{kilenc}
\end{eqnarray}
Observe the identity
\begin{equation}\label{tiz}
\sum_{u\in U} {1\over |A|}\sum_{i\in A}
\mu(i)\PP[\xi_i=u]= {1\over |A|} \sum_{i\in A}\mu(i)\sum_{u\in U}
\PP[\xi_i=u]
={1\over |A|}.
\end{equation}
Now identity (\ref{tiz}) implies  (\ref{tizenegy}) and
inequalities (\ref{nyolc}-\ref{kilenc}) imply inequality (\ref{weidiff}):
\begin{eqnarray} \label{tizenegy}
\min_{a\in A} r_a &= & {1\over |A|}+\sum_{u\in U}
 \max_{a\in A}\Biggl\{\mu(a)\PP[\xi_a=u] - {1\over |A|} \sum_{i\in A}
\mu(i)\PP[\xi_i=u]\Biggl\}\\
&\geq  & {1\over |A|} +{1\over |A|(|A|-1)}\sum_{u\in U}
\sum_{1\leq i < j \leq |A|}
\bigl|\mu(i)\PP[\xi_i=u]-\mu(j)\PP[\xi_j=u]\bigl|. \label{weidiff}
\end{eqnarray}
Fix an arbitrary $a,b\in A$, and set $Q=\sum_{u\in U} \bigl|\mu(a)
\PP[\xi_a=u]-\mu(b)\PP[\xi_b=u]\bigl|$. Define
\begin{eqnarray*}
U^+ & = & \biggl\{u\in U: \ \PP[\xi_a=u]>\PP[\xi_b=u]\biggl\},\\
U^= & = & \biggl\{u\in U: \ \PP[\xi_a=u]=\PP[\xi_b=u]\biggl\},\\
U^- & = & \biggl\{u\in U: \ \PP[\xi_a=u]<\PP[\xi_b=u]\biggl\}.
\end{eqnarray*}
Define further $A^+=\sum_{u\in U^+} \PP[\xi_a=u]$,
$A^-=\sum_{u\in U^-} \PP[\xi_a=u]$,

\noindent
$B^+=\sum_{u\in U^+}
\PP[\xi_b=u]$, $B^-=\sum_{u\in U^-}\PP[\xi_b=u]$. Observe that
\[
d(a,b)=\sum_{u\in U} \bigl|\PP[\xi_a=u]-\PP[\xi_b=u]\bigl|=A^+-B^++B^--A^-.
\]
On the other hand,
\[
A^+ + A^-= 1-\sum_{u\in U^=}\PP[\xi_a=u]=  1-\sum_{u\in
    U^=}\PP[\xi_b=u]=B^+ + B^-.
\]
  From the last two equations we conclude that $d(a,b)=2(A^+-B^+)=2(B^-
-A^-)$.
We finish the proof by setting a lower bound on $Q$ with a case analysis.
\begin{itemize}
\item
If $\mu(b)=\mu(a)$,  $Q=\mu(a)d(a,b)$.
\item If $\mu(b)>\mu(a)$,
$$Q\geq \mu(a)\sum_{u\in U^-} \PP[\xi_b=u]-\PP[\xi_a=u]={1\over 2} \mu(a)
d(a,b).$$
\item
If $\mu(b)<\mu(a)$,
$$Q\geq \mu(a)\sum_{u\in U^+} \PP[\xi_a=u] -\PP[\xi_b=u] =
 {1\over 2} \mu(a)d(a,b).$$
\end{itemize}

In all cases, we have $Q\geq {1\over 2} \mu(a)
d(a,b)$.
Returning to (\ref{weidiff}), we find
\begin{equation} \label{tizenot}
\sum_{1\leq i < j \leq |A|}\sum_{u\in U}
|\mu(i)\PP[\xi_i=u]-\mu(j)\PP[\xi_j=u]|\geq {1\over 2} \sum_{a\in A}
\mu(a)\sum_{b\in A} d(a,b),
\end{equation}
and through (\ref{tizenegy}), (\ref{weidiff}) and (\ref{tizenot}), we have 
\begin{eqnarray*}
\min_{a\in A} r_a 
& \geq & {1\over |A|} +{1\over 2|A|(|A|-1)}  \sum_{a\in A}
\mu(a)\sum_{b\in A} d(a,b)\\
& \geq & {1\over |A|} +{1\over 2|A|(|A|-1)} \min_{a\in A} \sum_{b\in A}
d(a,b).
\end{eqnarray*}
\end{proof}

\section{The parametric setting: Maximum Likelihood Estimation (MLE)}
\label{three}
In this section we reconsider the question of how many i.i.d. samples are
required in order for parametric maximum likelihood to accurately
recover elements of a finite set.

Assume $B=\{(a,\theta):a\in A,
\theta\in\Theta(a)\}$, and\/ $\Xi: B\rightarrow U$ is a parametric
random function, where $A$ and $U$ are
finite sets. Define 
\begin{eqnarray} \label{u}
U^{+}& := & \{u: \PP[\xi_{(a,\theta)} =u]>0\},\\
\alpha &:=& \alpha_{(a,\theta)} = \min_{u \in  \label{alph}
U^{+}}\{\PP[\xi_{(a,\theta)} =
u]\},
\end{eqnarray}
and assume
\begin{equation}
\label{dtheta}
d:= d_{(a,\theta)} = \inf_{b \neq a, \theta' \in \Theta(b)} \sum_{u \in
   U}|\PP[\xi_{(a,\theta)} =u] - \PP[\xi_{(b,\theta')}=u]|>0.
\end{equation}
In our earlier work, Theorem 5  in \cite{inv2},    we showed that for
\begin{equation}
\label{suffice}
k \geq  f(\alpha, d)\log(\frac{2|U^{+}|}{\epsilon}),
\end{equation}
$k$ samples suffice
to reconstruct $a \in A$, from $(a, \theta)$ with probability at least
$1-\epsilon$ using MLE, more formally, for $\Xi^{(k)}:B\rightarrow U^k$,
$[R^{(k)}]'_{(a, \theta)}\geq 1-\epsilon$.
Our function $f$ in (\ref{suffice}) tends to infinity when either (or
both) $\alpha \rightarrow 0$ or
$d \rightarrow 0$. This dependence on $d$ is
reasonable (though not always necessary, see Section~\ref{sup}), however
the dependence on $\alpha$ is not clear, and raises two questions.
\begin{itemize}
\item[Q1] Is there an bound on $k$ (like (\ref{suffice})) but which depends
only on $|U^{+}|, \epsilon$ and $d$ and not on
$\alpha$?

\item[Q2] Moreover, can the function $f$ in  (\ref{suffice}) be replaced by a
function of just 
  $d$ and $\epsilon$ (and not $\alpha$ and $U^+$)  so that the
  resulting function  is still a valid bound
for $k$?
\end{itemize}

In this section we show that the answer to the first question is `yes'
(Theorem~\ref{mainpara}) while the answer
to the second is `no' (Example~\ref{examp}).

We begin by introducing some further notation.
For any two probability distributions $p, p'$ on a set $U$ let
$d_{KL}(p,p') = \sum_{u \in U: p_u>0} p_u\log(\frac{p_u}{p'_u}) \in [0, \infty)
\cup \{\infty\}$ denote the Kullback-Leibler distance of $p$ and $p'$,
and recall the standard inequality:
\begin{equation}
d_{KL}(p,p') \geq \frac{1}{2}d(p, p')^2,
\label{ineqKL}
\end{equation}
where $d(p,p')$ denotes as usual the variational distance, $ \sum_{u
   \in U}|p_u - p'_u|$. We will also use $d_2(p,p')=\Bigl(\sum_{u
   \in U}|p_u - p'_u|^2\Bigl)^{1/2}.$

\begin{lemma}
\label{seqlem}
Let $X_1, X_2, \ldots, X_k$ be a sequence of i.i.d. random variables taking
values in a finite set $U$.
Assume further that if $X_i$ takes a value with probability zero, then it
never takes this value. 
For each $u \in U$, let $\hat{p}_u:= \frac{1}{k} \sum_{i=1}^k \II(X_i = 
u)$ (the
normalized multinomial counts) and let $p_u = \PP[X_1=u]$.
Let $U^{+}:=\{u: p_u >0\}$. Then,
\begin{itemize}
\item[(i)]
$\PP[d_{KL}(\hat{p}, p)\geq \delta] \leq \frac{|U^{+}|}{k\delta}$,
\item[(ii)]
$\PP[d(\hat{p}, p) \geq \delta] \leq \frac{|U^{+}|}{k\delta^2}$.
\end{itemize}
\end{lemma}

\begin{proof}
{\em Part (i)}
Let $\hat{\Delta}_u = \hat{p}_u -p_u$. For $u \in U^+$, set $\hat{Q}_u =0$ if 
$\hat{p}_u=0$, while if
$\hat{p}_u>0$ set
\begin{eqnarray}\hat{Q}_u  & : =&\hat{p}_u\log(\frac{\hat{p}_u}{p_u}) = 
(p_u+\hat{\Delta}_u)\log(1+ \nonumber
\frac{\hat{\Delta}_u}{p_u})\\
& \leq & (p_u+\hat{\Delta}_u) \cdot \frac{\hat{\Delta}_u}{p_u} \label{qu}
= \hat{\Delta}_u+ \frac{\hat{\Delta}_u^2}{p_u}.
\end{eqnarray}
  Recall Markov's inequality, which states that if $X$ is non-negative random
variable, and $a>0$, then
\begin{equation}\label{markov}
\PP[X \geq a] \leq \frac{\EE[X]}{a}.
\end{equation}
Note that $\EE[(\hat{p}_u - p_u)^2] 
=Var[\hat{p}_u]=
  \frac{p_u(1-p_u)}{k}$, and 
applying (\ref{markov}) to
$$X= \sum_{u \in U^{+}}\frac{\hat{\Delta}_u^2}{p_u} \geq 0$$ and noting that
$\EE[X] = \frac{|U^{+}|-1}{k}$ 
gives $\PP[X \geq \delta] \leq \frac{|U^{+}|}{k\delta}$.
By definition, $d_{KL}(\hat{p}, p) = \sum_{u:\hat{p}_u \not= 0} \hat{Q}_u =\sum_
{u \in U^{+} } \hat{Q}_u$
and this 
is less
or equal to $X$ (by (\ref{qu}), and the identity
$\sum_{u \in U^+} \hat{\Delta}_u =0$),
which leads to the required inequality.

{\em Part (ii)}
By the Cauchy-Schwartz inequality,
$d^2(\hat{p}, p) \leq d_2^2(\hat{p}, p)\cdot|U^{+}|$ and so,
$$\PP\bigl[d(\hat{p},p) \geq \delta\bigl] \leq \PP\bigl[d_2^2(\hat{p}, p) \geq
\delta^2/|U^{+}|\bigl]
\leq \frac{|U^{+}|}{\delta^2}\EE\bigl[d_2^2(\hat{p}, p)\bigl],$$
by Markov's inequality (\ref{markov}).
 Now, $$\EE[d_2^2(\hat{p}, p)] = \EE[\sum_{u \in
U}(\hat{p}_u -p_u)^2] = \sum_{u \in U} Var[\hat{p}_u] =
\sum_{u \in U} \frac{1}{k}p_u(1-p_u) \leq \frac{1}{k}.$$
\end{proof}

\begin{corollary}
\label{corostreet}
Under the assumptions of Lemma~\ref{seqlem},
if $\delta<1$, $\epsilon>0$  and $k \geq 
\frac{2|U^{+}|}{\epsilon\delta^2}$, then
with probability at least
$1-\epsilon$, the 
inequalities $d_{KL}(\hat{p}, p)< \delta$ and
 $d(\hat{p}, p) < \delta$ simultaneously hold.
\end{corollary}

\begin{theorem}
\label{mainpara} Assume $B=\{(a,\theta):a\in A,
\theta\in\Theta(a)\}$, and\/ $\Xi: B\rightarrow U$ is a parametric
random function, where $A$ and $U$ are
finite sets. Recall definition (\ref{u}) and condition (\ref{dtheta}).
Provided $k \geq \frac{c_1|U^{+}|}{\epsilon d_{(a,\theta)}^4}$ with
$c_1=\frac{2}{(2-\sqrt{3})^2}$, the probability that MLE correctly
returns $a$ from $\Xi^{(k)}$
is at least $1-\epsilon$, i.e. $[R^{(k)}]'_{(a, \theta)}\geq 1-\epsilon$.
\end{theorem}

\begin{proof}
Let $p$ be the probability distribution on $U$ induced by $\xi_{(a,\theta)}$,
$c=2-\sqrt{3}$,
$E$ be the event that $d(\hat{p}, p) \leq c \cdot d_{(a,\theta)}$. For the
probability
distribution $q$ induced by $\xi_{(b, \theta')}$ where $b \neq a$, 
by the triangle inequality we have
$$d(\hat{p}, q) \geq |d(p,q) - d(\hat{p}, p)|.$$
Now, by assumption $d(p,q) \geq d_{(a,\theta)}$, and so, conditional on $E$,
$d(\hat{p}, q) \geq (1-c)d_{(a,\theta)}$.
Invoking the inequality (\ref{ineqKL}) gives
$$d_{KL}(\hat{p}, q) \geq \frac{1}{2}d(\hat{p},q)^2 \geq
\frac{1}{2}(1-c)^2d_{(a,\theta)}^2.$$
Thus, conditional on $E$ we have:
\begin{equation}
\label{likelin}
\sum_{u \in U^+}\hat{p}_u\log q_u \leq 
\sum_{u \in U^+}\hat{p}_u\log \hat{p}_u -
\frac{1}{2}(1-c)^2d_{(a,\theta)}^2.
\end{equation}
For $ x\in A,\omega\in
\Theta(x)$, consider \begin{equation} L(x,\omega)=\sum_{u\in U^+} \hat
p(u) \log \PP[\xi_{x,\omega}=u].  \end{equation} 
$L(x,\omega)$ is ${1\over k} $ times the natural logarithm of
the probability that the observed sequence of $U$-elements came from
$(x,\omega)$. Therefore $L(x,\omega)\leq 0$ is proportional to the
log-likelihood of $(x,\omega)$. Now consider the log likelihood ratio
$$\Delta L:= L(a,\theta)-L(b,\theta')=
  \sum_{u \in U^+} \hat{p}_u\log(p_u/q_u).$$
Conditional on  $E$ we have, by (\ref{likelin}),
\begin{equation}
\label{likelq}
\Delta L \geq  - \sum_{u \in U^+} \hat{p}_u \log(\frac{\hat{p}_u}{p_u})
+ \frac{1}{2}(1-c)^2d_{(a,\theta)}^2=\frac{1}{2}(1-c)^2d_{(a,\theta)}^2
-d_{KL}(\hat{p},p).
\end{equation}
So if we select
$\delta = c \cdot d_{(a,\theta)}^2$ in Corollary~\ref{corostreet} we 
can ensure
that with probability at least
$1-\epsilon$ that event $E$ occurs and also (since  $\frac{1}{2}(1-c)^2 = c$)
that
$d_{KL}(\hat{p}, p) < \delta= c \cdot d_{(a,\theta)}^2=
\frac{1}{2}(1-c)^2  d_{(a,\theta)}^2$, and so, by
(\ref{likelq}) we have $\Delta L >0$.  The value of $k$ that
Corollary~\ref{corostreet} requires is precisely that given
in the statement of this theorem. This completes the proof.
\end{proof}

{\bf Remarks}
\begin{itemize}
\item
Theorem~\ref{mainpara} also implies that for MLE in the {\em non}--parametric
setting, the number $k$ of
i.i.d. samples required to reconstruct an element $a \in A$ correctly with
probability
at least $1-\epsilon$ is bounded above by a function that depends just on
$|U^{+}|, \epsilon$ and $d_a := \min_{b \neq a} d(a,b)$. In \cite{inv1} 
an upper
bound
on $k$ was also derived, however it depended just on
$|A|, \epsilon$ and $d_a$.  Comparing these results suggests an 
interesting question: Is
there an upper bound
for $k$ (in the non-parametric setting) which depends just on $d_a$ and
$\epsilon$?

\item
We show below that the linear dependence of $k$ on $|U^+|$ in
Theorem~\ref{mainpara} is best possible in the sense
that no sublinear dependence is possible.
It is possible however that the exponent of 4 for 
 $d$ in Theorem~\ref{mainpara} 
might be reduced.

\end{itemize}

\subsection{Construction to show that $k$ must grow 
 linearly with
   $|U^+|$}
\label{examp}

We now show that Theorem~\ref{mainpara} cannot be improved by replacing
the dependence of $k$ on $|U^{+}|$ with a sublinear function
(like the  logarithmic dependence on $|U|^+$ in Theorem 5.1 \cite{inv2}), 
even when $d_{(a,\theta)}$ and $\epsilon$
are held constant.

Let $A = \{a,b\}$, with $\Theta(a)=\{*\}$, and $$\Theta(b) =\{\theta =
(\lambda_1, \ldots, \lambda_n): \ 
\sum_{i=1}^n \lambda_i=1, \forall i\  \lambda_i \geq 0  \}.$$
 Let $U=\{0,1, \ldots, n\}$.
Fix $\delta>0$ and consider the random function $\Xi$ defined as follows.
\begin{align*}
\PP[\xi_{(a,*)} =u] = \begin{cases} \delta, & \mbox{if $u=0$;}
\\ \frac{1-\delta}{n}, &
\mbox{if $u \in \{1, \ldots, n\}$;} \end{cases} \end{align*}
\begin{align*}
\PP[\xi_{(b,(\lambda_1, \ldots, \lambda_n))} =u] = \begin{cases} 2\delta, &
\mbox{if $u=0$;}
\\ \lambda_u(1-2\delta), &
\mbox{ if $u  \in \{1,\ldots, n\}$.} \end{cases} \end{align*}
We assume that $k\leq n$, otherwise we have nothing to prove.
For ${\bf u} = (u_1, \ldots, u_k) \in U^k$, let $x({\bf u}) = |\{i \in \{1,
\ldots, k\}: u_i = 0\}|$.
We have:
$$L_1:= \sup_{\theta \in \Theta(a)} \PP[\xi^{(k)}_{(a,\theta)}={\bf u}]
= \delta^{x({\bf u})}\left(\frac{1-\delta}{n}\right)^{k-x({\bf u})},$$
and
\begin{equation} \label{szam}
L_2 := \sup_{\theta' \in \Theta(b)} \PP[\xi^{(k)}_{(b,\theta')}={\bf u}]
\geq (2\delta)^{x({\bf u})}\left(\frac{1-2\delta}{k-x({\bf 
u})}\right)^{k-x({\bf
u})},
\end{equation}
since we are free to select $\theta \in \Theta(b)$ to be the uniform
distribution on $\{1,\ldots, n\}$
for those $i$ for which $u_i \neq 0$.
We will select $\delta$
sufficient small that
\begin{equation}
2(1-2\delta)^{\delta/2} >1.
\label{epseq}
\end{equation}
Now, suppose we generate $u$ randomly from $(a,*)$. 
Note that the value of $d_{(a,*)}$ is at least $\delta$, since
$$d((a,*), (b, \theta')) \geq |\PP[\xi_{(a,*)} = 0]-\PP[\xi_{(b,\theta')}
 =0]| =
\delta.$$

Then MLE will (incorrectly) reconstruct $b$
whenever $R: = L_2/L_1 >1$. We will show that this occurs with probability at
least $1-\epsilon$,
if $k$ is less than  $\frac{1}{2}|U^{+}|$, for any $\delta$
satisfying (\ref{epseq}) and any sufficiently large  $|U^{+}|$.

Note that by replacing $L_2$ by its lower bound (\ref{szam}), we can write
$R \geq Y^k$ where
$$Y = 2^{\rho}\left[\frac{n}{k} \cdot
\frac{(1-2\delta)}{(1-\delta)(1-\rho)}\right]^{1-\rho},$$
where $\rho := x({\bf u})/k$.
Now, if $k \leq \frac{1}{2} n$, then  since $((1-\delta)(1-\rho))^{-(1-\rho)}
\geq 1$,
$$Y \geq 2(1-2\delta)^{1-\rho}.$$
Now, for $\delta, \epsilon$ fixed, there exists a value of $k$, for which, with
probability at least
$1-\epsilon$, we have $\rho > \frac{1}{2}\delta.$  Thus for this value of $k$,
and any  $n> 2k$
inequality (\ref{epseq}) gives
$$Y \geq 2(1-2\delta)^{\delta/2} >1,$$
and so $R >1$; that is MLE will make an incorrect decision. Thus, we 
must have
$k \geq \frac{1}{2}n = \frac{1}{2}(|U^+|-1)$
in order to avoid this.

\subsection{Example to show that parametric MLE can still succeed when
   variational distance vanishes on each element of $A$}
\label{sup}

In the {\em non}-parametric setting, given a random function $\Xi: A
\rightarrow U$,  suppose that $d(a,b) = 0$ for two
elements $a,b \in A$. Then for {\em any} random function $\Gamma: U
\rightarrow A$ it is easily shown (eg.  by Theorem 3.1 of \cite{inv2}) that
\begin{equation}
\label{halfbound}
\min\{\PP[\gamma_{\xi_{a_1}} = a_1], \PP[\gamma_{\xi_{a_2}}
= a_2]\} \leq \frac{1}{2}.
\end{equation}
That is, if  the probability
distribution induced by $a_1$ and $a_2$ is the same, no method can recover
both $a_1$ and $a_2$ more accurately than by a toss of a fair coin.
We can ask if a
similar result holds for {\em parametric} MLE. That is, suppose that $A=
\{a_1, a_2\}$ and for a value $\theta_1 \in \Theta(a_1)$, and $\theta_2
\in \Theta(a_2)$ we have
\begin{equation}
\label{vanish}
d_{(a_1,\theta_1)} = d_{(a_2, \theta_2)} = 0,
\end{equation}
where $d_{(a,\theta)}$ is defined as in (\ref{dtheta}).
Note that Theorem~\ref{mainpara} does not give a finite bound on $k$
for MLE to accurately reconstruct $a_1$ or $a_2$. However it turns out
that for certain random functions satisfying (\ref{vanish}), if parametric
MLE is used to estimate
$a_1$ and $a_2$ from $k$ independent trials, then for
 any parameter $(a_i, \theta_i)$ 
chosen, and for even $k$,  the probability that the selection is correct 
is always strictly greater than $\frac{1}{2}$, moreover in all but one 
choice of
 the parameter settings (for $a_1$)
the probability the selection is correct tends to $1$ as
 $k \rightarrow \infty$ 
(in the other setting it tends to $\frac{1}{2}$ from above). 
For this example th
ere is a more pedestrian approach for estimating $a_1$ or $a_2$ from
the $k$ independent trials, for which the probability of making the 
correct reconstruction tends to $1$ as $k$ tends to infinity, for all parameter
 settings (in
 contrast to MLE which has problems at one particular
 parameter settings -- this
 illustrates again the care required in consistency arguments for MLE).
  Note also that in this example, with any parameters
 $\theta_1, \theta_2$, $d\biggl((a_1, \theta_1),
(a_2, \theta_2)\biggl) >0$ holds.

Let $A=\{a_1, a_2\}$, $U=\{(1,0), (1,1), (2,0), (2,1)\}$,
$\Theta(a_1) = [ \pi/4, 3\pi/4) $,
and $\Theta(a_2) =  (\pi/4, 3\pi/4].$
For $t\in \Theta(a_1)$, let $\PP[\xi_{(a_1,t)}=
(1, \lfloor 2t/\pi\rfloor)]=\sin^2 t$,
$\PP[\xi_{(a_1,t)}=(2,\lfloor 2t/\pi\rfloor)]=\cos^2 t$;
and 
for $t\in \Theta(a_2)$, let $\PP[\xi_{(a_2,t)}=
(1, \lfloor 2t/\pi\rfloor)]=\cos^2 t$,
$\PP[\xi_{(a_2,t)}=(2,\lfloor 2t/\pi\rfloor)]=\sin^2 t$.

The key observation for the argument that follows is
 that $\sin^2 t > \cos^2 t$ 
in
$(\pi/4, 3\pi/4)$, while in the endpoints $\sin^2 t =1/2= \cos^2 t$.
It is easy to see that 
$\lim_{t\rightarrow {\frac{\pi}{4}}^+} d 
\biggl( (a_1, \pi/4), (a_2, t)\biggl)=0$, and hence
$d_{(a_1, \pi/4)}=0$. A similar argument shows that 
$d_{(a_2, 3\pi/4)}=0$. It is
 also easy to
see that the distributions of all $\xi_{(a_i,t)}$ random 
variables are different. The only
possible problem would be the distributions of $\xi_{(a_1,\pi/4)}$
and $\xi_{(a_2,3\pi/4)}$-- however in this case we have the
 second coordinates in the elements of $U$ to
separate these distributions. There is a 
 pedestrian way to guess where an element of $U$ came from. 
Count the ones and twos
in the first coordinates after $k$ independent trials. If there 
are more ones,
 then select $a_1$, if there are more $2$'s then select $a_2$, while 
in the case of a tie, 
 if $\lfloor 2t/\pi \rfloor=0$, then select $a_1$, otherwise 
select $a_2$.
(note that 
$\lfloor 2t/\pi \rfloor=0$ is constant over the trials). MLE pretty
 much does the same, the only
thing that requires more careful analysis is whether MLE correctly 
returns
$(a_1, \pi/4)$ and $(a_2, 3\pi/4)$. Focus on $(a_1, \pi/4)$, as
 the other problem is analogous.
Let \# 1 and \# 2 denote the number of  ones and twos
in the first coordinates in $\xi^{(k)}_{(a_1, \pi/4)   }$. Let $p$ be 
the probability of the
event  $X_1=$ ``\# 1 $>$ \# 2'';  by symmetry it is also the probability 
of the 
event 
$X_2=$ ``\# 1 $<$ \# 2'',
and let $q$ be the probability of the event
 $X_3=$ ``\# 1 $=$ \# 2''. Note that 
MLE
correctly returns $a_1$ for events $X_1$ and $X_3$ (but not for $X_2)$,
 and hence
$[R^{(k)}]'_{ (a_1, \pi/4)  }\geq p+q=\frac{1+q}{2}> \frac{1}{2}$.
The claim holds for $X_3$ for the following reason. The probability that 
$\xi^{(k)}_{(a_1, \pi/4)}$ yields the particular observed $k$-sequence
 conditional on $X_3$ is
$2^{-k}$, while the probability that 
$(a_2, \theta_2)$ generated the particular 
observed $k$-sequence conditional on event $X_3$ 
is $p^{k/2}(1-p)^{k/2}$
for some $p\not= 1/2$, and this second probability 
is strictly smaller than $2^{-k}$.

Informally, the reason for this
phenomena is that
the parameter space associated to $a_i$ is tuned for `fitting' data
that is produced by the pair $(a_i, \theta_i)$.

Despite this somewhat surprising result, one can easily derive a
parametric analogue of (\ref{halfbound}) for any random function $\Xi: B
\rightarrow U$ (where $B = \{(a, \theta):
\theta \in \Theta(a)\}$ as usual)
under the stronger condition
that $d((a_1, \theta_1), (a_2, \theta_2))=0$ where
$d((a_1, \theta_1), (a_2, \theta_2))$ is the variational distance between the
distributions of the
$U$--valued random variables
$\xi_{(a_1, \theta_1)}$ and $\xi_{(a_2, \theta_2)}$.  In this case, for any
random function (not just parametric MLE)
$\Gamma \rightarrow U$ that is independent of $\Xi$ it is easily shown that
$$\min\{\PP[\gamma_{\xi_{(a_1, \theta_1)}} = a_1], \PP[\gamma_{\xi_{(a_2, 
\theta_2)}}
= a_2]\} \leq \frac{1}{2}.$$
Of course this bound applies also for $k$ i.i.d. trial experiments.

\subsection{Application of Theorem~\ref{mainpara}}

As a simple illustration of the use of Theorem~\ref{mainpara}, we describe an
application to the reconstruction of
phylogenetic trees from binary sequences according to a simple Markov process
(the CFN model). Such processes are
central to much of molecular biology (see eg. \cite{fels}). Let $A$ denote the
three binary phylogenetic trees that have leaf set
$X=\{1,2,3,4\}$. For a tree $T=(V_T,E_T) \in A$, $\Theta(a)$ is the set of
functions $p:E_T \rightarrow [0,0.5]$ which
assign to each edge $e$ of $T$ an associated {\em substitution probability}.
Under the CFN model a state is assigned
uniformly at random to a leaf (eg. leaf 1) and states are assigned recursively
to the remaining vertices of the tree
by (independently) changing the state ($0$ to $1$ or $1$ to $0$) across each
edge $e$ of $T$ with probability $p(e)$. This gives
a (marginal) probability distribution on each of the 16 site patterns $c: X
\rightarrow \{0,1\}$ (further details concerning
this model can be found in \cite{inv2} or \cite{sem03}). Thus if we 
generate $k$
site patterns i.i.d. from
the pair $(T, p)$ we can ask how large $k$ must be in order for MLE to
accurately reconstruct
$T$.  To ensure that $d_{(T,p)}>0$ one must impose the following condition on
$p$.
\begin{itemize}
\item[(P)]
For each of the four edges $e$ of $T$ incident with a leaf we have
$p(e) \leq g< \frac{1}{2}$; and
for the central edge $e$ of $T$,  $p(e) \geq f >0$.
\end{itemize}

 From \cite{teasing} (Lemma 6.3) we have $d_{(T,p)} \geq H(f,g)>0$ for a
continuous function $H$.
Note that condition (P) can allow arbitrarily small values for
$\alpha_{(T,p)}:(= \min_{u \in U^{+}}\{\PP[\xi_{(T,p)} = u]\}$
even when $f$ and $g$ take fixed values (since condition (P) allows two 
adjacent
edges incident with leaves of $T$ to
both have arbitrarily small $p(e)$ values, and the probability of any site
pattern that assigns these two leaves
different states can therefore be made as close to zero as we wish).
Consequently, the main result from \cite{inv2} does not provide any (finite)
estimate for the site patterns
required for MLE to correctly reconstruct a tree.  However we may
applying Theorem~\ref{mainpara} in this setting, and since $|U^{+}| 
\leq 16$, we
obtain
an explicit upper bound on the
number of site patterns required to reconstruct each
phylogenetic tree on four leaves correctly with probability at least $1 -
\epsilon$.

\section{Acknowledgments}
We would like to thank Linyuan (Lincoln) Lu, for suggesting a shorter
proof of the
implication $(i) \Rightarrow (iii)$ of Proposition~\ref{main1}.


\begin{thebibliography}{99}

\bibitem{casellaberger} G. Casella and R.  L. Berger,  
{\em Statistical  Inference},
The Wadsworth \& Brooks/Cole Statistics/Probability Series, 
Wadsworth \& Brooks/Cole Advanced Books \& Software,
Pacific Grove, CA, 1990.


\bibitem{chang} J. T. Chang,  Full reconstruction of
   Markov models on evolutionary trees: identifiability and consistency,
   {\em Math. Biosci.} {\bf 137} (1996)  51--73.




\bibitem{everitt98}
B. S. Everitt,    {\em The Cambridge Dictionary of Statistics},
 Cambridge Univ. Press,
Cambridge, UK, 1998.




\bibitem{far99}
 J. S. Farris, 
  Likelihood and inconsistency, {\em Cladistics} {\bf 15} (1999)
 199--204. 


\bibitem{fels}  J. Felsenstein,  {\em Inferring Phylogenies,}
 Sinauer Press, 2004.


\bibitem{rog97}
J. S. Rogers,  On the consistency of maximum likelihood estimation 
of phylogenetic trees from nucleotide sequences,
 {\em Syst. Biol.} {\bf 46}
(1997) 354--357.

\bibitem{rog01}
J. S. Rogers,  Maximum likelihood estimation of phylogenetic
 trees is consistent when substitution rates vary according to the
 invariable sites plus gamma
distribution, {\em Syst. Biol.} {\bf 50} {2001} 713--722.


\bibitem{schrijver} A. Schrijver,
{\em Theory of Linear and Integer Programming},
Wiley-Interscience Series in Discrete Mathematics, John Wiley \& Sons
Ltd., Chichester, 1986.

\bibitem{sem03} C. Semple,  and M.  Steel,  
{\em Phylogenetics.} Oxford Univ. Press, 2003.


\bibitem{sid98}
M. E.  Siddall,
  Success of parsimony in the four-taxon case: long-branch
 repulsion by likelihood in the Farris zone,
 {\em Cladistics} {\bf 14} (1998) 209--220. 

\bibitem{inv1} M. A.   Steel and  L. A. Sz{\'e}kely, 
   Inverting random functions, {\em Annals of Combinatorics}, {\bf 3} (1999)
 103--113.


\bibitem{inv2} M. A.  Steel and  L. A. Sz{\'e}kely, 
 Inverting random functions II:
explicit bounds for the discrete maximum likelihood estimation, with
applications, {\em SIAM J. Discr. Math.}
{\bf 15(4)} (2002) 562--575.

\bibitem{teasing} M.A. Steel and L.A. and Sz{\'e}kely, Teasing 
apart two
trees.
(submitted). See IMI Technical Reports 05:08 
{\tt http://www.math.sc.edu/\~{}IMI/technical/tech05.html}, 2005.

\bibitem{wald}  A. Wald,
A note on the consistency of the maximum likelihood estimate, {\em Ann. Math.
Stat.},
{\bf 20} (1949)  595--600.



\bibitem{yan94}
Z. Yang,  Statistical properties of the maximum likelihood method of
 phylogenetic estimation and comparison with distance matrix methods,
 {\em Syst.
 Biol.} {\bf 43} (1994) 329--342.

\bibitem{yan96}
Z. Yang, Phylogenetic analysis using parsimony and likelihood methods,
 {\em Journal of Molecular Evolution} {\bf 42} (1996) 1641--1650.

\end{thebibliography}
\end{document}